\documentclass[letterpaper, 10 pt, draftcls, onecolumn]{ieeeconf} 
\IEEEoverridecommandlockouts 
\usepackage{graphicx}
\usepackage{subfigure}
\usepackage{amsmath}
\usepackage{amsfonts}
\usepackage[psamsfonts]{amssymb}
\usepackage{cases}
\usepackage{mathtools}
\usepackage{color}
\usepackage{epstopdf}
\usepackage{dsfont}
\usepackage{mathtools}
\usepackage{bbm}
\usepackage{cite}
\usepackage{diagbox}

\newtheorem{definition}{Definition}

\newtheorem{lemma}{Lemma}
\newtheorem{corollary}{Corollary}

\newtheorem{theorem}{Theorem}

\newtheorem{problem}{Problem}

\newcommand{\R}{\mathbb{R}}
\newcommand{\N}{\mathbb{N}}

\newcommand{\E}{\mathbb{E}}
\renewcommand{\P}{\mathbb{P}}

\newcommand{\nodeset}{V}
\newcommand{\nodenum}{n}
\newcommand{\edgeset}{E}
\newcommand{\graph}{G}
\newcommand{\randset}{\mathcal{G}}
\newcommand{\lapset}{\mathcal{L}}
\newcommand{\union}{U}
\newcommand{\fatone}{\mathbbm{1}}
\newcommand{\unionset}{\mathcal{U}}
\newcommand{\var}{\textnormal{Var}}
\newcommand{\eig}{\textnormal{eig}}

\mathtoolsset{showonlyrefs=true}

\makeatletter
\renewcommand{\fnum@figure}{Figure \thefigure}
\makeatother

\title{\LARGE \bf
On the Connectivity of Unions of Random Graphs
}

\author{Matthew T. Hale$^\star$
\thanks{$^\star$Department of Mechanical and Aerospace Engineering, University of Florida,
Gainesville, FL, USA. Email: \texttt{matthewhale@ufl.edu}.}
}

\begin{document}
\maketitle
\thispagestyle{empty}
\pagestyle{empty}

\begin{abstract}
Graph-theoretic tools and techniques have seen wide use in the multi-agent
systems literature, and the unpredictable nature of some multi-agent 
communications has been successfully modeled using random communication graphs.
Across both network control and network optimization, a common assumption is
that the union of agents' communication graphs is connected across any finite interval
of some prescribed length, and some convergence results explicitly depend
upon this length. Despite the prevalence of this assumption and the prevalence
of random graphs in studying multi-agent systems, 
to the best of our knowledge, there has not been a 
study dedicated to determining how many random graphs must
be in a union before it is connected. 
To address this point, this paper solves two related problems. 
The first bounds the number of random graphs required in a
union before its expected algebraic connectivity exceeds 
the minimum needed for connectedness. 
The second bounds the probability that a union of random graphs is connected.
The random graph model used is the Erd\H{o}s-R\'{e}nyi model, and,
in solving these problems, 
we also bound the expectation and
variance of the algebraic connectivity of unions of such graphs.
Numerical results for several use cases are given to supplement the theoretical
developments made. 
\end{abstract}

\section{Introduction}
Multi-agent systems have been studied in a number
of applications, including sensor networks \cite{cortes02}, 
robotics \cite{soltero13}, communications \cite{chiang07}, 
and smart power grids \cite{caron10}. 
Across these
applications, the agents in a network and
their associated communications are often
abstractly represented as graphs \cite{mesbahi10}. 
In general, graph-theoretic methods in multi-agent systems
represent each agent as a node in a graph and
each communication link as an edge, and multi-agent
coordination algorithms have been developed for both static
and time-varying graphs \cite[Chapter 1.4]{mesbahi10}. 

Time-varying random graphs in particular have been used to model communications
which are unreliable and intermittent due to interference
and poor channel quality \cite[Chapter 5]{mesbahi10}, and
such graphs have seen use in a number of multi-agent settings. 
For example, distributed agreement problems over random graphs are studied
in \cite{touri11} and \cite{hatano05},
while optimization over random graphs was explored in \cite{lobel11}.
The work in \cite{yazicioglu15} provides a means to modify random graphs
to make them robust to network failures, and \cite{lewis13} discusses
general properties of random graphs as they pertain to multi-agent systems.  
A broad
survey of graph-theoretic results for control can be found in \cite{mesbahi10},
and well-known graph-theoretic results in optimization
include \cite{blondel05,nedic09,zhu12}. 

When time-varying graphs (random or not) are used, a common assumption is that
the unions of these graphs are connected over intervals
of some finite length, i.e., the graph
containing all edges present over time is itself a
connected graph. A partial
sampling of works using this assumption (or a related variant)
includes 
\cite{moreau05,ren05,tseng90,tseng91,ren05b,chen12,kia15,blondel05,jadbabaie03,touri09,nedic10,nedic09,olshevsky11,zhu12,feyzmahdavian14,nedic07}. 
In addition, some works derive
convergence rates or other results that explicitly use the length of such
intervals, including \cite{feyzmahdavian14,tseng90,nedic07,tseng91,chen12,blondel05,touri09,olshevsky11,nedic09}. 
In applying these results, one may
wish to determine the time needed for the system to attain a connected union graph. 
To the best of our knowledge, no study
has been undertaken that addresses this problem for unions
of random graphs, despite their frequent use in multi-agent systems. 

Owing to the success of Erd\H{o}s-R\'{e}nyi graphs in
modeling some time-varying multi-agent communications \cite[Chapter 5]{mesbahi10}, 
we consider unions of random graphs generated
by the Erd\H{o}s-R\'{e}nyi model and examine the
connectedness of such unions. 
In particular, this paper solves two problems:
lower-bounding the number of random graphs required in a union before
one may expect it to be connected (in a precise sense
to be defined in Section~\ref{sec:review}), and
lower-bounding the probability that a union of random graphs
is connected.

Our results use spectral properties of the first four 
(matrix-valued) moments of the Laplacian of a union of random graphs.
The eigenvalues of these moments are used to bound the expected value of the 
Laplacian's second-smallest eigenvalue, called the 
\emph{algebraic connectivity} \cite{fiedler73}
of the underlying union graph. 
This bound in turn enables a lower bound on the number
of graphs needed in a union to have its algebraic connectivity
reach a specified expected value, and 
also enables a lower bound
on the probability of the algebraic connectivity exceeding some given threshold. 

The results presented rely heavily upon the spectral properties of 
random graphs' Laplacians, which are random matrices. 
A common approach to analyzing the spectra of random matrices
is to let the dimension of the matrix grow arbitrarily large
\cite{diaconis94,furedi81,wigner58,juhasz82}, and
the work in \cite{coja07} considers similar asymptotic results
focused on Laplacians of random graphs. 
For random graphs specifically, a common approach is to derive results
in which the size of the graph
grows arbitrarily large, and doing so
enables results that hold for \emph{almost all} graphs \cite{bollobas01}.
While there is clear theoretical
appeal to such results, 
our focus on multi-agent systems leads us to consider
non-asymptotic results precisely because such systems are
typically comprised by a fixed number of agents.
Our results are therefore stated for graphs of fixed
(but unspecified) size. 

In addition, while some work on
random graphs considers edge probabilities that bear some known
relationship to the number of nodes in a graph \cite{krivelevich03,tran13}, 
we do not do so here. Our use of random graphs to model multi-agent
communications is inspired by applications in which poor channel
quality, interference, and other factors make communications
unreliable. In such cases, the probability of a communication
link being active
may not bear any known relationship to the size of the network.
We therefore proceed with edge probabilities and network sizes
that are fixed and not assumed to be related. 

The rest of the paper is organized as follows. 
Section~\ref{sec:review} reviews the required elements
of graph theory and gives formal statements for the two
problems that are the focus of this paper.
Then, Section~\ref{sec:moments} computes moments of
random graph Laplacians and certain spectral properties
of these moments to enable the results of
Section~\ref{sec:algconn}. Section~\ref{sec:algconn}
then presents the main results of the paper and solves
the problems stated in Section~\ref{sec:review}. 
Next, Section~\ref{sec:simulation} presents numerical
solutions to several instantiations of the problems
studied. Finally, Section~\ref{sec:conclusion}
provides concluding remarks and future directions
for extending this work.

\section{Review of Graph Theory and Problem Statements} \label{sec:review}
In this section, we review the required elements of 
graph theory. We begin with basic definitions, including
the definition of algebraic connectivity, and
then review the Erd\H{o}s-R\'{e}nyi model for
random graphs; throughout this paper, all uses of the phrase 
``random graphs'' refer to Erd\H{o}s-R\'{e}nyi graphs. 
Then we formally state the two problems
solved in this paper. Below, we use the notation
$[m] := \{1, \ldots, m\}$ for any $m \in \N$. 

\subsection{Basic Graph Theory}
A \emph{graph} is defined over a set of 
nodes, denoted $\nodeset$, and describes connections
between these nodes in the form of edges, which
are contained in an edge set $\edgeset$. Formally,
for $\nodenum$ nodes, $\nodenum \in \N$, the elements of $\nodeset$
are indexed over $[\nodenum]$. The set of edges
in the graph is a subset
\begin{equation}
\edgeset \subseteq \nodeset \times \nodeset,
\end{equation}
where a pair $(i, j) \in \edgeset$ if nodes $i$ and $j$
share a connection, and $(i, j) \not\in \edgeset$ if they do not. 
This paper considers graphs which are \emph{undirected}, meaning an 
edge $(i, j)$ is not distinguished from an edge $(j, i)$, 
and \emph{simple}, so that $(i, i) \not\in \edgeset$ for all
$i$. A graph $\graph$ is then defined as the pair
$\graph = (\nodeset, \edgeset)$. One main focus of this paper
is on \emph{connected} graphs, which we define now. 

\begin{definition} (E.g., \cite{godsil01})
A graph $\graph$ is called \emph{connected} if, for all $i \in [\nodenum]$ and $j \in [\nodenum]$,
$i \neq j$, there is a sequence of edges one can traverse from node $i$ to node $j$, i.e.,
there is a sequence of indices $\{i_{\ell}\}_{\ell=1}^{k}$ and nodes
$\{v_{p}\}_{p=i_1}^{i_k}$ such that $\edgeset$ contains all of the edges
\begin{equation}
(i, v_{i_1}), \: (v_{i_1}, v_{i_2}), \: (v_{i_2}, v_{i_3}), \: \ldots, \: (v_{i_{k-1}}, v_{i_k}), \: (v_{i_k}, j). 
\end{equation}
\hfill $\triangle$
\end{definition}

The results of this paper are developed in terms of
graph Laplacians, which we define now. First, the adjacency
matrix $A(\graph) \in \R^{\nodenum \times \nodenum}$ associated with the graph $\graph$ is defined
element-wise as
\begin{equation}
a_{ij} = 
\begin{cases}
1 & (i, j) \in \edgeset \\
0 & \textnormal{otherwise},
\end{cases}
\end{equation}
where $a_{ij}$ is the $i^{th}j^{th}$ element of $A(\graph)$. 
When there is no ambiguity, we will simply denote $A(\graph)$ by $A$.
Because we consider undirected graphs, $A$ is symmetric by definition. 

Next, the degree matrix $D(\graph) \in \R^{\nodenum \times \nodenum}$ associated with a graph 
$\graph$ is a diagonal matrix
whose entries count the number of edges connecting to a node.
Using $d_i$ to denote the \emph{degree} of node $i$, we find
\begin{equation}
d_i = \sum_{\substack{j=1 \\ j \neq i}}^{\nodenum} a_{ij} = |\{j \mid (i, j) \in \edgeset\}|,
\end{equation}
where $|\cdot|$ denotes the cardinality of a set. Then the degree matrix
associated with a graph $\graph$ is 
\begin{equation}
D(\graph) = \textnormal{diag}(d_1, d_2, \ldots, d_n)
\end{equation}
which we will denote $D$ when $\graph$ is clear from context.
Clearly $D$ is also symmetric by definition. 

The Laplacian of a graph $\graph$ is then defined as
\begin{equation}
L(\graph) = D(\graph) - A(\graph),
\end{equation}
which will be written simply as $L$ when $\graph$ is unambiguous. 
The results of this paper rely in particular on spectral properties
of $L$. Letting $\lambda_k(\cdot)$ denote the $k^{th}$ smallest
eigenvalue of a matrix, it is known that
$\lambda_1(L) = 0$ for all graph Laplacians \cite{mesbahi10}, and thus
we have
\begin{equation} \label{eq:lambdai}
0 = \lambda_1(L) \leq \lambda_2(L) \leq \cdots \leq \lambda_{\nodenum}(L).
\end{equation}
The value of $\lambda_2(L)$ is central to the work in this paper and
some other works in graph theory, and it gives rise to the following
definition. 

\begin{definition} \label{def:lambda2} (From \cite{fiedler73})
The \emph{algebraic connectivity} of a graph $\graph$ is the second
smallest eigenvalue of its Laplacian, $\lambda_2(L)$,
and $\graph$ is connected if and only if $\lambda_2(L) > 0$.
\hfill $\triangle$
\end{definition}

This paper is dedicated to studying the statistical properties of $\lambda_2$
for unions of random graphs. Toward doing so, we now review the 
necessary elements of the theory of random graphs. 

\subsection{Random Graphs}
Several well-known random graph models exist in the literature \cite{erdos59,watts98},
and Erd\H{o}s-R\'{e}nyi graphs in particular have been successfully used in the multi-agent systems
literature. Erd\H{o}s-R\'{e}nyi graphs can model, for example, unreliable, intermittent
and time-varying communications in multi-agent networks \cite{mesbahi10},
and we therefore consider the Erd\H{o}s-R\'{e}nyi model in this paper. Under this model, a graph on 
$\nodenum$ vertices contains each admissible edge with some fixed \emph{edge probability} $p \in (0, 1)$.
Therefore, for each $i \in [\nodenum]$ and $j \in [\nodenum]$ with
$i \neq j$, an Erd\H{o}s-R\'{e}nyi graph satisfies
\begin{equation} \label{eq:pijE}
\P[(i, j) \in \edgeset] = p \textnormal{ and } \P[(i, j) \not\in E] = 1 - p.
\end{equation}
Equivalently, based on Equation~\eqref{eq:pijE} one finds
\begin{equation}
\E[a_{ij} = 1] = p \textnormal{ and } \E[a_{ij} = 0] = 1 - p,
\end{equation}
i.e., that $a_{ij}$ is a Bernoulli random variable for $i \neq j$. 

We denote the sample space of all Erd\H{o}s-R\'{e}nyi graphs on $\nodenum$ nodes with edge probability $p$
by $\randset(n, p)$, and we denote the set of Laplacians of all such graphs 
by $\lapset(n, p)$. One approach to spectral graph theory commonly used in the literature is
to let $\nodenum \to \infty$ \cite{bollobas98}. The value of doing so is that one may draw conclusions that
hold for \emph{almost all} graphs in a rigorous way. In the study of multi-agent systems, however,
one is often focused on networks with a fixed number of agents that is not well approximated
by letting $\nodenum$ become arbitrarily large. Accordingly, 
we develop our results in terms of an arbitrary but fixed value of $\nodenum$. 

In addition,
some well-known results in the graph theory literature assume that $p$ has some
known relationship to $\nodenum$ \cite{bollobas98}, or else that the number of edges
in a random graph has some relationship to the number of nodes in the graph \cite{erdos60}. While the theoretical
utility of these relationships is certainly clear from those works, this relationship will often not hold in 
multi-agent systems where communications are unpredictable because these communications
are affected by a wide variety of external factors. We therefore proceed with a value of $p \in (0, 1)$ that
is not assumed to have any relationship to the value of $\nodenum$. 

In the study of multi-agent systems, it is also common for algorithms and results to be stated in terms
of unions of graphs, which we define now. 
\begin{definition}
For a collection of graphs ${\{\graph_k = (V, E_k)\}_{k=1}^{N}}$ defined on the same node set $\nodeset$,
the union of these graphs, denoted $\union_{N}$, is defined as
\begin{equation}
\union_{N} := \bigcup_{k=1}^{N} G_k = \left(V, \cup_{k=1}^{N} E_k\right),
\end{equation}
i.e., the \emph{union graph} $\union_N$ contains all edges in all $N$ graphs
that comprise the union. \hfill $\triangle$ 
\end{definition}

\subsection{Problem Statement}
A common requirement in some multi-agent systems is that the communication graphs in a network
form a connected union graph over intervals of some fixed length. To help determine
when this occurs, we formulate and solve two related problems in this paper.
The first concerns when
a union graph has expected algebraic connectivity above some threshold. 
\begin{problem} \label{prob:exp}
Find $N \in \N$ such that $\E\big[\lambda_2\big(\union_N\big)\big] \geq \lambda_{min}$, 
where $\lambda_{min}$ is the minimum algebraic connectivity of a connected
graph, $\union_N$ is given by
\begin{equation}
\union_N = \bigcup_{k=1}^{N} \graph_k,
\end{equation}
and $\graph_k \in \randset(n, p)$ for all $k \in [N]$. \hfill $\lozenge$
\end{problem}

The second problem we solve concerns the probability with which a union graph has algebraic
connectivity exceeding the minimum among connected graphs. 
\begin{problem} \label{prob:prob}
Given $N \in \N$, lower bound the value of 
\begin{equation}
\P\left[\lambda_2\big(\union_{N}) \geq \lambda_{min}\right],
\end{equation}
where $\lambda_{min}$ is the minimum algebraic connectivity of a connected
graph, and 
where $\union_{N}$ is defined as it is in Problem~\ref{prob:exp}. \hfill $\lozenge$
\end{problem}
 
Section~\ref{sec:moments} next provides theoretical developments that enable
the solutions to these problems in Section~\ref{sec:algconn}.

\section{Moments and Spectra of Random Graph Laplacians} \label{sec:moments}
Towards solving Problems~\ref{prob:exp} and \ref{prob:prob}, 
this section computes the first four moments of a random
graph's Laplacian. The values of these moments will be used
below to compute the expectation and variance of
the algebraic connectivity of random graphs and their
unions, and these results later enable solutions to
Problems~\ref{prob:exp} and \ref{prob:prob}. 

\subsection{Moments of Random Graph Laplacians}
We begin by stating a lemma that will be used below to 
compute eigenvalues of moments of $L$. 

\begin{lemma} \label{lem:abeig}
Let $I$ be the $\nodenum \times \nodenum$ identity matrix,
and let $J$ by the $\nodenum \times \nodenum$ matrix 
whose entries are all $1$. Then the matrix
\begin{equation}
M := (\alpha - \beta)I + \beta J = \left(\begin{array}{ccccc} \alpha & \beta  & \beta  & \cdots & \beta \\
                                                               \beta & \alpha & \beta  & \cdots & \beta \\
                                                               \beta & \beta  & \alpha & \cdots & \beta \\
                                                              \vdots & \vdots & \vdots & \ddots & \vdots \\
                                                               \beta & \beta   & \beta  & \cdots & \alpha\end{array}\right)
\end{equation}
has $\alpha + (n-1)\beta$ as an eigenvalue with multiplicity one and $\alpha - \beta$ as an eigenvalue
with multiplicity $n-1$.
\end{lemma}
\emph{Proof:} See Lemma~1 in \cite{hale17}. \hfill $\blacksquare$

We now present the first four moments of a random graph Laplacian $L \in \lapset(n, p)$. 

\begin{lemma} \label{lem:moments}
Let $\graph \in \randset(n,p)$ have Laplacian $L$. Then
\begin{align}
\E[L]   &= p(nI - J) \\
\E[L^2] &= \big[(n-2)p^2 + 2p\big](nI - J) \\
\E[L^3] &= \big[(n-2)(n-4)p^3 + 6(n-2)p^2 + 4p\big](nI - J) \\
\E[L^4] &= \big[(n-7)(n-3)(n-2)p^4 + 6(2n-7)(n-2)p^3 + 25(n-2)p^2 + 8p\big](nI - J),
\end{align}
where $I$ is the $\nodenum \times \nodenum$ identity matrix and $J$ is the
$\nodenum \times \nodenum$ matrix of ones. 
\end{lemma}
\emph{Proof:} See Lemma~2 in \cite{hale17}. \hfill $\blacksquare$

Next, we present a lemma showing the equivalence between the expected spectrum of
powers of $L$ and the spectrum of the corresponding moments of $L$. The purpose of this lemma
is to enable the use of Lemma~\ref{lem:moments} in computing the expected
eigenvalues for a random graph Laplacian ${L \in \lapset(n, p)}$. 

Before doing so, we draw an important distinction between
the eigenvalues $\{\ell_i\}_{i=1}^{n}$ studied in this 
section and the eigenvalues $\{\lambda_i\}_{i=1}^{n}$
in Section~\ref{sec:review}. The eigenvalues
$\{\ell_i\}_{i=1}^{n}$ comprise an unordered collection
and are simply the eigenvalues
of a random graph's Laplacian; as a result, each
$\ell_i$ is itself a random variable. In the setting
of random graphs, $\lambda_i$ is then the $i^{th}$ \emph{order statistic}
over these random variables, i.e., the $i^{th}$ smallest
value realized by any of the random variables
in the collection $\{\ell_i\}_{i=1}^{n}$. 

More concretely,
using the convention that $\ell_1 = 0$ is always the zero
eigenvalue of a random graph's Laplacian
(which is guaranteed to exist by Equation~\eqref{eq:lambdai}), the algebraic
connectivity of a random graph is defined as
\begin{equation}
\lambda_2 = \min_{2 \leq i \leq \nodenum} \ell_i. 
\end{equation}
This section characterizes each $\ell_i$, and Section~\ref{sec:algconn} uses
these results to characterize $\lambda_2$. Towards doing so, we have
the following lemma. 

\begin{lemma} \label{lem:chi}
Let $L$ be the Laplacian of a random graph ${\graph \in \randset(n, p)}$, and
let $\eig(M)$ denote the set of eigenvalues of a matrix $M \in \R^{\nodenum \times \nodenum}$. Then,
for $k \in [4]$,
\begin{equation}
\eig\big(\E[L^k]\big) = \E\big[\eig(L^k)\big].
\end{equation}
\end{lemma}
\emph{Proof:} 
The expected value of any diagonal element
of a graph Laplacian takes the form
\begin{equation}
\E[L_{ii}] = \E\left[\sum_{\substack{j=1 \\ j \neq i}}^{n} a_{ij} \right] = (n-1)p
\end{equation}
because the random variables $a_{ij}$ are independent Bernoulli random variables
which take value $1$ with probability $p$. Summing these diagonal entries, we find the
expected trace of $L$ to be
\begin{equation} \label{eq:etrace}
\E[\textnormal{trace}(L)] = \sum_{i=1}^{n} \E[L_{ii}] = n(n-1)p.
\end{equation}

Denote the eigenvalues of a matrix $L \in \lapset(\nodenum, p)$
by $\ell_i$, $1 \leq i \leq n$. 
Due to the fact that all off-diagonal
elements of $L$ are i.i.d. random variables, and
that the diagonal
elements are simply sums of these variables, the non-zero eigenvalues
of $L$ are equal in expectation. That is, apart from the guaranteed 
zero eigenvalue $\ell_1 = 0$, all other eigenvalues
have equal expectation, precisely because all off-diagonal entries
of $L$ take the same form and because all diagonal entries
do as well. 
By Equation~\eqref{eq:etrace} and the 
definition of the trace of a matrix, we then find
\begin{equation}
\sum_{i=2}^{n} \E[\ell_i] = n(n-1)p,
\end{equation}
giving
\begin{equation} \label{eq:privacy_np1}
\E[\ell_i] = np.
\end{equation}

As for $\eig(\E[L])$, we note that
\begin{equation}
\E[L] = p(nI - J),
\end{equation}
which by Lemma~\ref{lem:abeig} has eigenvalues
\begin{equation} \label{eq:privacy_np2}
\ell_1 = 0 \textnormal{ and } \ell_i = np \textnormal{ for } k \in \{2, \ldots, n\}.
\end{equation}
Comparing Equations~\eqref{eq:privacy_np1} and \eqref{eq:privacy_np2}, we find that
\begin{equation}
\E[\eig(L)] = \eig(\E[L]).
\end{equation}

By the same reasoning, one can repeatedly exploit the symmetries
of $L$ and its powers to obtain the same result for $L^k$. 
\hfill $\blacksquare$

In words, Lemma~\ref{lem:chi} says that the expected spectrum of $L^k$
is equal to the spectrum of the expectation of $L^k$. It was shown 
in Lemma~\ref{lem:moments} that
$\E[L^k]$ takes a simple form for ${k \in [4]}$, and thus Lemma~\ref{lem:chi}
simplifies the process of computing the expected eigenvalues
of $L^k$. 

Using Lemmas~\ref{lem:moments} and \ref{lem:chi}, we can compute moments
of the eigenvalues of $L$. The eigenvalues of $L$ are
denoted by $\ell_i$, and these are random variables because they are functions
of the entries of $L$, which in turn are either Bernoulli random variables
(for off-diagonal entries) or sums of Bernoulli random variables (for diagonal
entries). We begin with the smallest eigenvalue of $L$. 

\begin{lemma} \label{lem:ell1}
Let $L \in \lapset(n, p)$. Then $\ell_1 = 0$ is
in $\E[\eig(L)]$, and $\ell_1^k = \ell_1 = 0$
is in $\E[\eig(L^k)]$ for $k \in \{2, 3, 4\}$. 
\end{lemma}
\emph{Proof:} Let $\fatone$ denote the vector in $\R^{\nodenum}$ whose
entries are all $1$. Then we see that 
\begin{equation}
(nI - J)\fatone = 0. 
.
\end{equation}
%
%
Then $0$ is an eigenvalue of $nI - J$.
Because $\E[L]$ in Lemma~\ref{lem:moments} is a scalar multiple
of $nI - J$, 
$0$ is an eigenvalue of $\E[L]$. By Lemma~\ref{lem:chi},
$0$ is then also the expected value of an eigenvalue of $L$. 
Because $0$ is an eigenvalue of $nI - J$, it is also
an eigenvalue of $(nI - J)^k$ for $k \in \{2, 3, 4\}$, and
repeating the preceding argument for these values of $k$
completes the lemma. 
\hfill $\blacksquare$
%

Having established the expectation of $\ell_1$ in Lemma~\ref{lem:ell1},
we now compute the first four moments of all other $\ell_i$'s. 

\begin{theorem} \label{thm:ellimoments}
Let $\graph \in \randset(\nodenum, p)$ and let $\ell_i$ denote the $i^{th}$
eigenvalue of its Laplacian. 
For all $i \in [\nodenum] \backslash \{1\}$, 
\begin{align}
\E[\ell_i]   &= np \\
\E[\ell_i^2] &= n(n-2)p^2 + 2np \\
\E[\ell_i^3] &= n(n-2)(n-4)p^3 + 6n(n-2)p^2 + 4np \\
\E[\ell_i^4] &= n(n-7)(n-3)(n-2)p^4 + 6n(2n-7)(n-2)p^3 + 25n(n-2)p^2 + 8np.
\end{align}
\end{theorem}
\emph{Proof:} Applying Lemma~\ref{lem:abeig} to Lemma~\ref{lem:moments}
gives the above quantities as eigenvalues of $\E[L]$, and Lemma~\ref{lem:chi}
establishes that these eigenvalues are moments of eigenvalues of $L$. 
\hfill $\blacksquare$

Using Theorem~\ref{thm:ellimoments}, we have the following
corollary which computes the variances of $\ell_i$ and $\ell_i^2$,
and these variances will be applied in the next section to
bound certain properties of $\lambda_2$.

\begin{corollary} \label{cor:variances}
Let $\graph \in \randset(\nodenum, p)$ and let $\ell_i$ denote
the $i^{th}$ eigenvalue of its Laplacian. Then for $i \in [\nodenum] \backslash \{1\}$,
\begin{equation}
\var[\ell_i] = 2npq
\end{equation}
and
\begin{equation}
\var[\ell_i^2] = n(n-7)(n-3)(n-2)p^4 + 6n(2n-7)(n-2)p^3 + 25n(n-2)p^2 + 8np - \big(n(n-2)p^2 + 2np\big)^2.
\end{equation}
\end{corollary}
\emph{Proof:} 
By definition,
\begin{equation}
\var[\ell_i] = \E[\ell_i^2] - \E[\ell_i]^2 \textnormal{ and } \var[\ell_i^2] = \E[\ell_i^4] - \E[\ell_i^2]^2,
\end{equation}
and the result follows using the results of Theorem~\ref{thm:ellimoments}. 
\hfill $\blacksquare$

Having characterized certain statistical properties of
the collection $\{\ell_i\}_{i=1}^{\nodenum}$, the next
section translates these properties into bounds
on statistical properties of $\lambda_2$.

\section{Algebraic Connectivity of Unions of Random Graphs} \label{sec:algconn}
This section translates the bounds on $\ell_i$ derived in Section~\ref{sec:moments}
for single random graphs into bounds on $\lambda_2$ for unions of random
graphs. First, we show that a union of random graphs can itself be represented
as a random graph with a different edge probability. Second, we present
results that bound the expectation of order statistics in terms of the
expectations and variances of the underlying collection of random variables.
Third, we present our solutions to Problems~\ref{prob:exp} and \ref{prob:prob}.
In this section, we use the notation $q = 1 - p$. 

\subsection{Unions of Random Graphs are Random Graphs}
Section~\ref{sec:moments} derived results for single random graphs, and we show
now that these results are easily adapted to unions of random graphs because
such unions are themselves equivalent to single random graphs with
a different edge probability.

\begin{lemma} \label{lem:unioneq}
Let $\unionset_N(n,p)$ denote the set of all unions of $N$ random graphs on $\nodenum$ 
nodes with edge probability $p$, i.e.,
\begin{equation}
\unionset_N(n, p) := \left\{\bigcup_{i=1}^{N} \graph_i \mid \graph_i \in \randset(n, p)\right\}. 
\end{equation}
Then
\begin{equation}
\unionset_N(n, p) = \randset\big(n, 1 - (1 - p)^N\big).
\end{equation}
\end{lemma}
\emph{Proof:} 
Consider some $G \in \mathcal{U}_{N}(n, p)$. Fix any admissible node indices $i$ and $j$. 
Then an edge is absent between $i$ and $j$ only if it is absent in all $N$ graphs that
comprise $G$. That is, an edge between $i$ and $j$ is absent in $G$ with probability
$q^N$. Then that edge is present with probability $1 - q^N = 1 - (1 - p)^N$.
\hfill $\blacksquare$

With Lemma~\ref{lem:unioneq}, results pertaining to individual random graphs can be
applied to unions of such graphs with only minor modifications. 

\subsection{Expectation of Order Statistics}
It was noted in Section~\ref{sec:moments} that the algebraic connectivity of a
random graph is the first order statistic over the non-zero eigenvalues of
that random graph's Laplacian, i.e.,
\begin{equation} \label{eq:l2ordstat}
\lambda_2 = \min_{2 \leq i \leq \nodenum} \ell_i.
\end{equation} 
Thus, while the expected value of each $\ell_i$ is known, the expected 
value of $\lambda_2$ is not. To apply what is known about $\ell_i$
to $\lambda_2$, we state the following lemma from \cite{arnold79}
which bounds the expectation of order statistics in terms of
properties of the underlying collection of random variables.

\begin{lemma} \label{lem:ordstat} 
Let $X_1$, $X_2$, $\ldots$, $X_m$ be jointly distributed with 
common mean $\mu$ and variance $\sigma^2$. Then the $k^{th}$ order statistic
of this collection, denoted $X_{k:m}$, has expectation bounded according to
\begin{equation}
\mu - \sigma\sqrt{\frac{m-k}{k}} \leq \E[X_{k:m}] \leq \mu + \sigma\sqrt{\frac{k-1}{n-k+1}}.
\end{equation}
\end{lemma}
\emph{Proof:} See Equation~(4) in \cite{arnold79}. 
\hfill $\blacksquare$

The bound in Lemma~\ref{lem:ordstat} is shown in \cite{arnold79}
to be tight when the underlying
random variables have identical means and variances.
Using Lemma~\ref{lem:ordstat}, we now bound the expected value of $\lambda_2$
for a single random graph.

\begin{lemma} \label{lem:l2graph}
Let $\graph \in \randset(\nodenum, p)$. Its algebraic connectivity, $\lambda_2$, has expectation
bounded according to
\begin{equation} \label{eq:leml2graphmain}
\max\{np - \sqrt{2n(n-2)pq}, 0\} \leq \E[\lambda_2] \leq np,
\end{equation}
where $q := 1 - p$. 
\end{lemma}
\emph{Proof:} This follows from Equation~\eqref{eq:l2ordstat} and
using Lemma~\ref{lem:ordstat} with $m = \nodenum-1$, $\mu$ from
Theorem~\ref{thm:ellimoments}, and $\var[\ell_i]$ from 
Corollary~\ref{cor:variances}. 
The non-negativity of the left-hand side of Equation~\eqref{eq:leml2graphmain}
follows from the non-negativity of all eigenvalues of all $L \in \lapset(\nodenum, p)$,
stated in Equation~\eqref{eq:lambdai}. 
\hfill $\blacksquare$

It is possible that the left-hand side of Equation~\eqref{eq:leml2graphmain}
is zero for some values of $p$. In particular, a straightforward calculation shows
that the left-hand side of Equation~\eqref{eq:leml2graphmain} is only positive
when
\begin{equation}
p > \frac{2n - 4}{3n - 4}
\end{equation}
and for $p$ outside this range, Lemma~\ref{lem:l2graph} does not provide a lower
bound on $\lambda_2$ beyond its non-negativity (which can be inferred from 
the non-negativity of each $\ell_i$). 
However, despite this limitation, Lemma~\ref{lem:l2graph} will be instrumental in solving 
Problem~\ref{prob:exp} below. Before doing so, we now
bound the variance of $\lambda_2$ 
by following an argument similar to that in Lemma~\ref{lem:l2graph}.

\begin{lemma} \label{lem:l2vargraph}
Let $\graph \in \randset(\nodenum, p)$. Its algebraic connectivity, $\lambda_2$,
has variance bounded according to
\begin{align}
\var[\lambda_2] &\leq n(n-2)p^2 + 2np - \big(np - \sqrt{2n(n-2)pq}\big)^2 \\
\var[\lambda_2] &\geq n(n-2)p^2 + 2np - \sigma[\ell_i^2]\sqrt{n-2} - n^2p^2,
\end{align}
where
\begin{equation}
\sigma[\ell_i^2] = \Big(n(n-7)(n-3)(n-2)p^4 + 6n(2n-7)(n-2)p^3 + 25n(n-2)p^2 + 8np - \big(n(n-2)p^2 + 2np\big)^2\Big)^{1/2}
\end{equation}
as in Corollary~\ref{cor:variances}. 
\end{lemma}
\emph{Proof:} 
Using Lemma~\ref{lem:ordstat} and Theorem~\ref{thm:ellimoments}, we find that
\begin{equation}
n(n-2)p^2 + 2np - \sigma[\ell_i^2]\sqrt{n-2} \leq \E[\lambda_2^2] \leq n(n-2)p^2 + 2np.
\end{equation}
Using that $\var[\lambda_2] = \E[\lambda_2^2] - \E[\lambda_2]^2$ and
the bounds on $\E[\lambda_2]$ from Lemma~\ref{lem:l2graph}, the result follows. 
\hfill $\blacksquare$

To assess connectivity of random graphs, the final result needed is a lower
bound on the algebraic connectivity of connected graphs. It is known that the
connected graph with least algebraic connectivity is a line graph \cite{fiedler73},
and we present this value below.

\begin{lemma} \label{lem:line}
The minimum algebraic connectivity attained by a connected graph on $\nodenum$ nodes is
that of a line graph, equal to
\begin{equation}
\lambda_{min} = 2\left(1 - \cos\frac{\pi}{n}\right).
\end{equation}
\end{lemma}
\emph{Proof:} See Proposition 1.12 in \cite{belhaiza05}. \hfill $\blacksquare$

While Definition~\ref{def:lambda2} says that a graph is connected if and only if $\lambda_2 > 0$,
Lemma~\ref{lem:line} shows that there is a minimum value of $\lambda_2$ attained by any 
connected graph. Definition~\ref{def:lambda2} still holds because any graph with $\lambda_2 > 0$
will also have $\lambda_2 \geq \lambda_{min}$. However, when computing $\E[\lambda_2]$ for a
random graph, it is possible to have $0 < \E[\lambda_2] < \lambda_{min}$, in which case
$\E[\lambda_2]$ is not large enough to imply connectivity of the underlying graph.
Thus, while an actual graph has $\lambda_2 \geq \lambda_{min}$ whenever $\lambda_2 > 0$,
an ``expected graph'' may not. Accordingly, our solutions to Problems~\ref{prob:exp}
and \ref{prob:prob} use $\lambda_{min}$ as the desired lower bound on $\lambda_2$, 
with the knowledge that doing so is sufficient for connectivity. 

\subsection{Solutions to Problems~\ref{prob:exp} and \ref{prob:prob}}
We now present the main results of the paper: solutions to Problems~\ref{prob:exp} and \ref{prob:prob}. 
We begin with Problem~\ref{prob:exp} and provide a lower bound on the number of random graphs
needed in a union before its expected algebraic connectivity is bounded below by the minimum
among all connected graphs, as determined in Lemma~\ref{lem:line}.

\begin{theorem} \label{thm:prob1} \emph{(Solution to Problem~\ref{prob:exp})}
The expected algebraic connectivity of a union of $N$ random graphs is bounded below by
the minimum for connected graphs if
\begin{equation}
N \geq N_{min} := \frac{1}{\log q}\log\left(\frac{4n^2 + 4n\cos\frac{\pi}{n} - \tau(n) - 8n}{6n^2 - 8n}\right),
\end{equation}
where
\begin{equation}
\tau(n) := \bigg(16n^2(n - 2)\left(1 - \cos\frac{\pi}{n}\right) + 32n(2 - n)\left(1 - \cos\frac{\pi}{n}\right)^2 + 4n^2(n-2)^2\bigg)^{1/2}.
\end{equation}
\end{theorem}
\emph{Proof:} From Lemma~\ref{lem:line}, the minimum algebraic connectivity attained by any connected graph is 
$\lambda_{min} = 2\left(1 - \cos\frac{\pi}{n}\right)$. 
By Lemma~\ref{lem:l2graph}, we find that
\begin{equation} \label{eq:main11}
np - \sqrt{2n(n-2)pq} \leq \E[\lambda_2],
\end{equation}
where $\lambda_2$ is the algebraic connectivity for a union of $N$ random graphs.

Using Lemma~\ref{lem:unioneq} and replacing $p$ by $\hat{p} := 1 - (1 - p)^N$ in Equation~\eqref{eq:main11}
gives
\begin{equation}
n\hat{p} - \sqrt{2n(n-2)\hat{p}\hat{q}} \leq \E[\lambda_2],
\end{equation}
where $\hat{q} = 1 - \hat{p}$. To lower-bound
$\E[\lambda_2]$ by $\lambda_{min}$, it is sufficient for
\begin{equation}
\lambda_{min} \leq n\hat{p} - \sqrt{2n(n-2)\hat{p}\hat{q}},
\end{equation}
or, rearranging terms, it is sufficient for
\begin{equation}
2n(n-2)\hat{p}\hat{q} \leq n^2\hat{p}^2 - 2n\hat{p}\lambda_{line} + \lambda_{line}^2. 
\end{equation}
Replacing $\hat{p}$ by $1 - \hat{q}$ gives a quadratic inequality for $\hat{q}$, which
can be solved for $\hat{q}$ using the quadratic equation. Then, expanding $\hat{q} = (1 - p)^N$
and solving for $N$ gives the desired result. 
\hfill $\blacksquare$

Section~\ref{sec:simulation} gives numerical lower bounds on $N$ generated by Theorem~\ref{thm:prob1}
for a range of values of $\nodenum$ and $p$. We emphasize that Theorem~\ref{thm:prob1} holds for
any fixed values of $\nodenum$ and $p$, without requiring any relationship between them.
In Theorem~\ref{thm:prob1}, it can also be shown that $\tau(n)$ is dominated by $2n^2$ for large $n$,
and therefore the argument of the second $\log$
is dominated by the $2n^2$ term in the numerator and the $6n^2$ term in the denominator, resulting
in this term limiting to $\log(1/3) = -\log(3)$. Therefore, as $n$ becomes large, the lower bound on $N_{min}$
in Theorem~\ref{thm:prob1} approaches a limiting value, namely, for large $n$,
\begin{equation} \label{eq:23sat}
N_{min} \approx -\frac{\log(3)}{\log q}.
\end{equation}

Having solved Problem~\ref{prob:exp}, we now focus on Problem~\ref{prob:prob}. 
Toward solving Problem~\ref{prob:prob}, we now state the Paley-Zygmund inequality in the form
in which we use it below. 

\begin{lemma} \label{lem:pzineq} \emph{(Paley-Zygmund inequality)}
Let $Z$ be a non-negative random variable with $\var[Z] < \infty$ and let $\theta \in [0, 1]$. 
Then
\begin{equation}
\P(z > \theta\E[Z]) \geq (1 - \theta)^2\frac{\E[Z]^2}{\E[Z^2]}.
\end{equation}
\end{lemma}
\emph{Proof:} See \cite{paley32}. \hfill $\blacksquare$

\begin{theorem} \label{thm:prob2} \emph{(Solution to Problem~\ref{prob:prob})}
The probability that the algebraic connectivity of a union of $N \geq N_{min}$ random graphs is 
at least the minimum algebraic connectivity of a connected graph is
\begin{equation}
\P[\lambda_2(U_N) \geq \lambda_{min}] \geq \left(1 - \frac{2\left(1 - \cos\frac{\pi}{n}\right)}{n\hat{p}}\right)^2\left(\frac{\left(n\hat{p} - \sqrt{2n(n-2)\hat{p}\hat{q}}\right)^2}{n(n-2)\hat{p}^2 + 2n\hat{p}}\right),
\end{equation}
where $\hat{p} = 1 - (1 - p)^N$ and $\hat{q} = 1 - \hat{p}$. 
\end{theorem}
\emph{Proof:} 
From Lemma~\ref{lem:l2graph}, we see that
\begin{equation}
n\hat{p} - \sqrt{2n(n-2)\hat{p}\hat{q}} \leq \E[\lambda_2],
\end{equation}
while from Lemma~\ref{lem:l2vargraph} we find that
\begin{equation}
\E[\lambda_2^2] \leq n(n-2)\hat{p}^2 + 2n\hat{p}. 
\end{equation}
Substituting these bounds into Lemma~\ref{lem:pzineq} and setting
${\theta = \lambda_{min}/\E[\lambda_2]}$
completes the proof. \hfill $\blacksquare$

The condition that $N \geq N_{min}$ in Theorem~\ref{thm:prob2} is enforced so that 
$\theta = \lambda_{min}/\E[\lambda_2]$ can be used in the Paley-Zygmund inequality.
Of course, this condition can be eliminated and a different form of probabilistic
bound can be derived in place of Theorem~\ref{thm:prob2}. As stated,
Theorem~\ref{thm:prob2} gives a probabilistic bound that is a function only
of $N$, the number of random graphs in a union, because $n$ and $p$ are fixed.
Together, Theorems~\ref{thm:prob1} and \ref{thm:prob2} characterize any union
of Erd\H{o}s-R\'{e}nyi graphs and help determine the number of 
graphs required to attain connectivity in such unions. 

In the next section, we give numerical results derived from Theorems~\ref{thm:prob1}
and \ref{thm:prob2} for values of $n$ and $p$ across several orders
of magnitude.

\section{Numerical Results} \label{sec:simulation}
In this section, we simulate the results of Theorems~\ref{thm:prob1}
and \ref{thm:prob2} to provide numerical solutions to Problems~\ref{prob:exp}
and \ref{prob:prob} for select values of $n$ and $p$. 

\subsection{Numerical Results for Problem~\ref{prob:exp}}
We now present results using Theorem~\ref{thm:prob1} by 
providing values of $N_{min}$ for a range of values of $n$ and
$p$, representing the solutions to Problem~\ref{prob:exp}
under different conditions. 

Table~\ref{tab:sol1} gives the value of $N_{min}$ (rounded
up) as determined by Theorem~\ref{thm:prob1} for $n$ ranging
from $10$ to $100,000$ and $p$ ranging from $0.00001$ to $0.1$. 
The values of $N_{min}$ shown in the table are the numbers of
graphs from $\randset(\nodenum, p)$ needed in a union
before the algebraic connectivity of the union has expectation
bounded below by that of a line graph (which has least
algebraic connectivity among all connected graphs). 

\begin{table}[htbp]
\centering
\begin{tabular}{|c|c|c|c|c|c|} \hline
\diagbox{$p$}{$n$} & $10$      & $100$     & $1,000$   & $10,000$  & $100,000$ \\ \hline\hline
$0.00001$          & $117,846$ & $110,539$ & $109,928$ & $109,868$ & $109,862$ \\ \hline
$0.0001$           & $11,785$  & $11,054$  & $10,993$  & $10,987$  & $10,986$  \\ \hline
$0.001$            & $1,178$   & $1,105$   & $1,099$   & $1,099$   & $1.099$   \\ \hline
$0.01$             & $118$     & $110$     & $110$     & $110$     & $110$      \\ \hline
$0.1$              & $12$      & $11$      & $11$      & $11$      & $11$       \\ \hline
\end{tabular}
\caption{Values of $N_{min}$, as determined by Theorem~\ref{thm:prob1}}
\label{tab:sol1}
\end{table}

Throughout these values, it can be seen that, for each fixed value of $n$,
an order of magnitude increase
in $p$ corresponds well to an order of magnitude decrease in the number
of graphs required for connectivity of a union. 
In addition, for a fixed value of $p$, increasing $n$ causes a decrease
in the lower bound on $N_{min}$. This means that, as graphs become larger,
fewer total graphs are required in a union to make it connected.
This occurs because a graph on $n$ nodes has $\frac{n(n-1)}{2}$ possible
edges and, because a larger graph has more possible edges, larger
graphs have more possible ways to attain connectivity, resulting
in fewer required in a union to make it connected. The limiting
behavior of Theorem~\ref{thm:prob1} seen in Equation~\eqref{eq:23sat} 
can also be seen in Table~\ref{tab:sol1},
where the lower bounds on $N_{min}$ appear to saturate when $p$
is held fixed and $n$ is increased. For example, for $p = 0.00001$, we have
\begin{equation}
-\frac{\log(3)}{\log q} = 109,861,
\end{equation}
which agrees closely with the values of $N_{min}$ for ${p = 0.00001}$ and
$n \geq 1,000$ seen in Table~\ref{tab:sol1}. 

\subsection{Numerical Results for Problem~\ref{prob:prob}}
We now present numerical results from Theorem~\ref{thm:prob2}. 
In particular, for $n = 50$ nodes and 
${p \in \{0.05, 0.10, 0.15, 0.20, 0.25\}}$, we present lower
bounds on the probability of a union of $N = 50$ random graphs
being connected. While Theorem~\ref{thm:prob1} concerns the 
expected algebraic connectivity of a union of random graphs,
Theorem~\ref{thm:prob2} concerns the probability
of the algebraic connectivity itself exceeding
that of a line graph (which is the minimum among all connected
graphs). Below, Table~\ref{tab:sol2} gives
numerical values for the probability of a particular union
of $N = 50$ random graphs having algebraic connectivity bounded
below by that of a line graph. 

\begin{table}[htbp]
\centering
\begin{tabular}{|c|c|c|c|c|c|} \hline
\diagbox{$n$}{$p$} & $0.05$  & $0.10$  & $0.15$  & $0.20$  & $0.25$ \\ \hline\hline
$50$               & $0.359$ & $0.810$ & $0.953$ & $0.989$ & $0.998$ \\ \hline
\end{tabular}
\caption{A lower bound on the probability of $N = 50$ graphs from $\randset(50, p)$ being connected,
as determined by Theorem~\ref{thm:prob2}}
\label{tab:sol2}
\end{table}

In Table~\ref{tab:sol2} we see that the probability of having a connected union increases
rapidly as $p$ increases, and that this probability approaches $1$ even when $p$ is 
far from $1$. 
Next, in Table~\ref{tab:sol3}, we present results that
have fixed values of $n=50$ and $p=0.10$, but changing values of $N$.

\begin{table}[htbp]
\centering
\begin{tabular}{|c|c|c|c|c|c|} \hline
\diagbox{$n$}{$N$} & $25$    & $50$    & $75$    & $100$   & $125$ \\ \hline\hline
$50$               & $0.377$ & $0.810$ & $0.947$ & $0.986$ & $0.996$ \\ \hline
\end{tabular}
\caption{A lower bound on the probability of $N$ graphs from $\randset(50, 0.1)$ being connected,
as determined by Theorem~\ref{thm:prob2}}
\label{tab:sol3}
\end{table}

Similar to what was seen in Table~\ref{tab:sol2}, these results show that the probability
of a union being connected increases rapidly with $N$. 
To further illustrate this trend, a union of $N = 250$ graphs with $n = 50$ and $p = 0.1$
is connected with probability at least $0.9998$ according to Theorem~\ref{thm:prob2}.

The numerical results in this section indicate that the results in 
Theorems~\ref{thm:prob1} and \ref{thm:prob2} solve Problems~\ref{prob:exp}
and \ref{prob:prob} in a manner which readily provides numerical results.
The three parameters $n$, $p$, and $N$ can vary dramatically across problem
formulations, though the results obtained here apply to broad
ranges across all three of these parameters, 
allowing for these results to be used in a wide range of applications.

\section{Conclusion} \label{sec:conclusion}
We presented results that determine the number of random
graphs required for their union to attain some lower
bound on its expected algebraic connectivity, and 
results that
lower-bound the probability with which a union of 
random graphs is connected. 
In multi-agent systems, a common assumption is that
agents' communication graphs have connected unions
over time, and these results can be used to enforce
this assumption.

Future work includes extending
our results to the cases of different edge probabilities
and of time-varying edge probabilities. Another direction
for future work concerns reformulating these results
for directed graphs, including with time-varying
probabilities that are direction-dependent.
The directed case breaks the symmetries used in 
this paper and would therefore likely require 
different techniques to derive bounds of the
same form derived here. 

\bibliographystyle{plain}{}
\bibliography{sources}

\end{document}